\documentclass[11pt]{article}
\usepackage[latin1]{inputenc}
\usepackage{epsfig}
\usepackage{color}
\usepackage[british,english]{babel}
\usepackage{amsthm}
\usepackage{amsmath}
\usepackage{amsfonts}
\usepackage{amssymb}
\usepackage{graphicx}
\newtheorem{corollary}{Corollary}[section]
\newtheorem{definition}[corollary]{Definition}

\newtheorem{lemma}[corollary]{Lemma}
\newtheorem{proposition}[corollary]{Proposition}
\newtheorem{remark}[corollary]{Remark}
\newtheorem{theorem}[corollary]{Theorem}
\newfont{\sBlackboard}{msbm10 scaled 900}

\newcommand{\mylabel}[1]{\label{#1}
            \ifx\undefined\stillediting
            \else \fbox{$#1$}\fi }
\newcommand{\BE}{\begin{equation}}
\newcommand{\BEQ}[1]{\BE\mylabel{#1}}
\newcommand{\EEQ}{\end{equation}}
\newcommand{\rfb}[1]{\mbox{\rm
   (\ref{#1})}\ifx\undefined\stillediting\else:\fbox{$#1$}\fi}
\newcommand{\half}   {{\frac{1}{2}}}

\newfont{\Blackboard}{msbm10 scaled 1200}
\newcommand{\bl}[1]{\mbox{\Blackboard #1}}
\newfont{\roma}{cmr10 scaled 1200}
\def\RR{{\rm I~\hspace{-1.15ex}R} }
\def\CC{\rm \hbox{C\kern-.56em\raise.4ex
         \hbox{$\scriptscriptstyle |$}\kern+0.5 em }}
\newcommand{\be}{\begin{equation}}
\newcommand{\ee}{\end{equation}}
\newcommand{\beq}{\begin{eqnarray}}
\newcommand{\eeq}{\end{eqnarray}}
\newcommand{\beqs}{\begin{eqnarray*}}
\newcommand{\eeqs}{\end{eqnarray*}}
\newcommand{\bt}{\begin{Theorem}}
\newcommand{\et}{\end{Theorem}}
\newcommand{\br}{\begin{remark}}
\newcommand{\er}{\end{remark}}
\newcommand{\bc}{\begin{Corollary}}
\newcommand{\ec}{\end{Corollary}}
\newcommand{\el}{\end{Lemma}}
\newcommand{\bd}{\begin{definition}}
\newcommand{\ed}{\end{definition}}

\newcommand{\nline}  {{\bl N}}
\newcommand{\rline}  {{\bl R}}

\def\cD{{\cal D}}

\newcommand{\Dscr} {{\cal D}}

\newcommand{\mm}    {{\hbox{\hskip 0.5pt}}}
\newcommand{\m}     {{\hbox{\hskip 1pt}}}
\newcommand{\nm}    {{\hbox{\hskip -3pt}}}

\newcommand{\bluff} {{\hbox{\raise 15pt \hbox{\mm}}}}

\newcommand{\FORALL} {{\hbox{$\hskip 11mm \forall \;$}}}

\newcommand{\rarrow} {{\,\rightarrow\,}}

\makeatletter
\def\section{\@startsection {section}{1}{\z@}{-3.5ex plus -1ex minus
    -.2ex}{2.3ex plus .2ex}{\large\bf}}
\makeatother
\def\be{\begin{equation}}
\def\ee{\end{equation}}

\def\ds{\displaystyle}
\begin{document}

\thispagestyle{empty}
\title{\bf Stabilization by switching control methods}
\author{Ka\"{\i}s AMMARI
\thanks{D\'epartement de Math\'ematiques,
Facult\'e des Sciences de Monastir, Universit\'e de Monastir, 5019 Monastir, Tunisie,
e-mail : kais.ammari@fsm.rnu.tn}
\, , Serge NICAISE \thanks{Universit\'e de Valenciennes et du Hainaut Cambr\'esis, LAMAV, FR CNRS 2956,  59313 Valenciennes Cedex 9, France, e-mail: snicaise@univ-valenciennes.fr}
\, and \, Cristina PIGNOTTI \thanks{Dipartimento di Matematica Pura e Applicata, Universit\`a di L'Aquila, Via Vetoio, Loc. Coppito, 67010 L'Aquila, Italy, \,
e-mail : pignotti@univaq.it}}
\date{}
\maketitle
{\bf Abstract.} {\small
In this paper we consider some stabilization problems for the wave equation with switching.
We prove  exponential stability results for appropriate damping coefficients. The proof of the main results is based on D'Alembert formula and some energy estimates.} \\


\section{Introduction}
\setcounter{equation}{0}
Our main goal is to study the pointwise or boundary stabilization of a switching delay wave equation in $(0,\ell)$. More precisely, we consider the
systems  given by :
\begin{align}
{}&u_{tt}(x,t)- u_{xx}(x,t)   = 0 ,  && \mbox{\rm in} \quad(0,\ell) \times (0,2 \ell),\label{a1}\\
{}&u_{tt}(x,t)- u_{xx}(x,t) + a \, u_t(\xi,t-2 \ell) \, \delta_\xi = 0 ,  && \mbox{\rm in } (0,\ell) \times (2 \ell,+\infty),\label{a1b}\\
{}&u(0,t)=0, \, u_x(\ell,t) = 0,&& \mbox{\rm on}\quad(0,+\infty), \label{a2}\\
{}&u(x,0) = u_0(x),  \quad  u_t(x,0)=u_1(x), && \mbox{\rm in} \quad (0,\ell),\label{a4}
\end{align}
and
\begin{eqnarray}
& &u_{tt}(x,t) - u_{xx}(x,t)=0\quad \mbox{\rm in}\quad(0,\ell) \times
(0,+\infty),\label{1.1}\\
& & u(0,t)  =0\quad \mbox{\rm on}\quad
(0,+\infty),\quad  \label{1.2}\\
& & u_x(\ell,t) = 0  \quad \mbox{\rm
on}\quad
(0,2\ell),\label{1.3}\\
& & u_x(\ell,t)= \mu_1 u_t(\ell,t)  \quad \mbox{\rm
on}\quad
(2(2i+1)\ell,2(2i+2)\ell), \forall i\in \nline,\label{dampstand}\\
& & u_x(\ell,t)= \mu_2 u_t(\ell, t-2\ell)  \quad
\mbox{\rm on}\quad
(2(2i+2)\ell,2(2i+3)\ell), \forall i\in \nline,\label{dampdelay}\\
& &u(x,0)=u_0(x)\quad \mbox{\rm and}\quad u_t(x,0)=u_1(x)\quad
\hbox{\rm in} \quad(0,\ell),\label{1.4}
\end{eqnarray}
where $\ell > 0, \mu_1,\mu_2, a$ and $\xi \in (0,\ell)$ are constants.
Here and below we denote by $\nline$ the set of the natural numbers while 
$\nline^*=\nline \setminus\{ 0\}.$

Note that in both cases, the feedbacks are unbounded.
 
Delay effects arise in many applications and practical problems and
it is well-known that an arbitrarily small delay may destabilize a
system which is uniformly asymptotically stable in absence of delay
(see e.g. \cite{Datko,DLP,Datko97},  \cite{NPSicon})). Nevertheless recent papers reveal that particular choice of delays may restitute exponential stability property, see \cite{g,gt,guo}.

We refer also to \cite{amman,ANP,NPSicon,NVCOCV10} for stability results for systems with time delay
due to the presence of ``good'' feedbacks compensating the destabilizing delay effect.  

Note that the above systems are exponentially stable in absence of
time delay, and if $\mu_1=\mu_2<0$ (see e.g.
\cite{Chen}) for the second system and if $a > 0,\  \xi$ admits a coprime factorization $\frac{p}{q}$ and $p$ is odd (the best rate is obtained for $\xi = \frac{\ell}{2}$, see e.g. \cite{aht}) for the first system.

In this paper we propose a new approach that consists to stabilize the wave system
by a control law that uses informations from the past (by  switching or not).
This means that the  stabilization is obtained  by a control method (that we propose to call 
switching control method)
and not by a feedback law.
For the first system this law is given by the term
$a \, u_t(\xi,t - 2\ell) \, \delta_\xi$ in 
 (\ref{a1b}) for $t\geq 2\ell$, while 
for the second system it corresponds to the term   $\mu_2 u_t(\ell,t-2 \ell)$) in a switched control form.
Using D'Alembert formula and some energy estimates,
we will show that for any $a\in (0,2)$   and $\xi = \frac{\ell}{2}$, system \eqref{a1}--\eqref{a4} 
is exponentially stable. On the other hand we show that appropriate choices of $\mu_1$
and $\mu_2$ yield the exponential  stability of  \rfb{1.1}--\rfb{1.4}.

The same approach is briefly treated in higher dimension,
here our approach combines observability estimates and   some energy estimates.

For the existence results, let us recall the following facts.
Let $A = - \partial^2_x$ be the unbounded operator in $H = L^2(0,\ell)$
with domain
$$H_{1}={\cal D}(A) = \left\{u \in H^2(0,\ell); \, u(0) = 0, \, u_x(\ell) = 0 \right\},$$
$$
H_\half = {\cal D}(A^\half) = \left\{u \in H^1(0,\ell); \, u(0) =0 \right\}.$$
We define 
$$
B_1 \in {\cal L}(\rline, H_{-\half}), \, B_1k = k \, \sqrt{a} \, \delta_\xi, \forall \, k \in \rline, \, B_1^*u = \sqrt{a} \, u(\xi), \, \forall \, u \in H_\half,
$$
and
$$
B_2 \in {\cal L}(\rline, H_{-\half}), \, B_2k = \sqrt{\mu_1} \, A_{-1} D k = k \sqrt{\mu_1} \delta_\ell, \forall \, k \in \rline, \, B^*u = \sqrt{\mu_1} \, u(\ell), \, \forall \, u \in H_\half,
$$
where $A_{-1}$ is the extension of $A$ to $H_{-1} = ({\cal D}(A))^\prime$ and $D$ is the Dirichlet map ($Dk = kx$ on $(0,\ell)$) and $H_{-\half} = (H_\half)^\prime$ (the duality is in the sense of $H$).

To study the well--posedness of the systems
\eqref{a1}--\eqref{a4} and \rfb{1.1}--\rfb{1.4}, we write them as an abstract
Cauchy problem in a product space, and use the semigroup approach.
For this purpose, take the Hilbert space ${\mathcal H}~:=~ H_{\half}
\times H$ and the unbounded linear
operators \be \label{op} {\mathcal  A} : {\mathcal D}({\mathcal
A}) =  H_{1} \times
H_{\half} \subset {\mathcal  H} \longrightarrow {\mathcal  H}, \,
{\mathcal A} \left(
\begin{array}{ccc}u_1\\u_2\end{array}\right) = \left(
\begin{array}{l}
u_2 \\
- Au_1 
\end{array}
\right), 
\ee
and
 $$ {\mathcal  A}_d : {\mathcal D}({\mathcal
A}_d) =  \left\{(u,v) \in [H^2(0,\ell) \times H^1 (0,\ell)] \cap {\cal H}\, : \, u_x(\ell) = \mu_2 \, v(\ell) \right\} \subset {\mathcal  H} \longrightarrow {\mathcal  H}, $$
\be \label{opb}
{\mathcal A}_d \left(
\begin{array}{ccc}u_1\\u_2\end{array}\right) = \left(
\begin{array}{l}
u_2 \\
- Au_1 
\end{array}
\right). 
\ee
It is well known that the operators $({\mathcal  A}, {\mathcal  D}({\mathcal  A}) )$ and $({\mathcal  A}_d, {\mathcal  D}({\mathcal  A}_d) )$ defined by \eqref{op} and (\ref{opb}), 
generate a strongly continuous semigroup of contractions on ${\mathcal H}$ denoted respectively  $({\mathcal
T}(t))_{t\geq0}$ (we also denote   $({\mathcal
T}_{-1}(t))_{t\geq0}$ the extension of $({\mathcal
T}(t))_{t\geq0}$ to $H_{-1}$), $({\mathcal
T}_d(t))_{t\geq0}$.
\begin{proposition}\label{propexistunic}
\begin{enumerate}
\item
The system
\rfb{a1}--\rfb{a4} is well-posed. More precisely, for every
$(u_0,u_1)\in {\mathcal  H}$, the solution of \rfb{a1}--\rfb{a4} is given by 
$$
\left(
\begin{array}{ccc}  u(t) \\ u_t(t)
\end{array}\right)  = \left\{
\begin{array}{ll}
\left(
\begin{array}{ccc}
u^0(t)\\
u_t^0(t)
\end{array}\right)
  = {\mathcal  T}(t)\left(
\begin{array}{ccc}u_0\\u_1\end{array}\right), \, 0 \leq t \leq 2 \ell, \\
\left(
\begin{array}{ccc}
u^{j}(t)\\
u_t^{j}(t)
\end{array}\right) = {\mathcal  T}(t-2j \ell)\left(
\begin{array}{ccc}u^{j-1}(2j \ell)\\u_t^{j-1}(2j \ell)\end{array}\right) + \\ \ds \int_{2j\ell}^{t} {\mathcal  T}_{-1}(t-s)\left(
\begin{array}{ccc}0\\- a \, u_t^{j-1}(s - 2 \ell) \delta_\xi\end{array}\right) \, ds, \, \\
2j \ell \leq t \leq 2 (j+1) \ell, j \geq 1.
\end{array}
\right.
$$ 
and satisfies $(u^j,u^j_t) \in C([2j\ell, 2(j+1)\ell], {\mathcal  H}), \, j \in \nline.$
\item
The system
\rfb{1.1}--\rfb{1.4} is well-posed. More precisely, for every
$(u_0,u_1)\in {\mathcal  H}$, the solution of \rfb{1.1}--\rfb{1.4} is given by 
$$
\left(
\begin{array}{ccc}
u(t)\\
u_t(t)
\end{array}\right)= \left\{
\begin{array}{ll}
\left(
\begin{array}{ccc}
u^0(t)\\
u_t^0(t)
\end{array}\right) = {\mathcal  T}(t)\left(
\begin{array}{ccc}u_0\\u_1\end{array}\right), \, 0 \leq t \leq 2 \ell, \\
\left(
\begin{array}{ccc}
u^{2j+1}(t)\\ u_t^{2j+1}(t)
\end{array}\right)
 = {\mathcal  T}_d(t-2 (2j+1) \ell)\left(
\begin{array}{ccc}u^{2j}(2 (2j+1) \ell )\\ u_t^{2j}(2 (2j+1) \ell )\end{array}\right), \\
\hspace{5cm} 2 (2j+1) \ell \leq t \leq 2 (2j +2)\ell, j \in \nline,\\
\left(
\begin{array}{ccc}
u^{2j+2}(t)\\ u_t^{2j+2}(t)
\end{array}\right) = {\mathcal  T}(t-2(2j+2) \ell)\left(
\begin{array}{ccc}u^{2j+1}(2(2j+2) \ell)\\u_t^{2j+1}(2(2j+2) \ell)\end{array}\right) + \\
\ds \int_{2(2j+2) \ell}^{t} {\mathcal  T}_{-1}(t-s)\left(
\begin{array}{ccc}0\\- \mu_1 \, u_t^{2j+1}(s- 2 \ell) \delta_\ell\end{array}\right) \, ds, \, \\
\hspace{5cm} 2(2j+2) \ell \leq t \leq 2 (2j+3) \ell, j \in \nline
\end{array}
\right.
$$ 
and satisfies $$(u^0,u^0_t) \in C([0, 2 \ell], {\mathcal  H}), \, (u^{2j+1},u^{2j+1}_t) \in C([2(2j+1)\ell, 2(2j+2)\ell], {\mathcal  H}), \,j \in \nline, \, $$ 
$$ (u^{2j+2},u^{2j+2}_t) \in C([2(2j+2)\ell, 2(2j+3)\ell], {\mathcal  H}), \, j \in \nline.$$
\end{enumerate}
\end{proposition}

For any solution of problem \rfb{a1}--\rfb{a4} respectively of  \rfb{1.1}--\rfb{1.4} we define the energy
\begin{equation}\label{energy}
\hspace{3cm}
E_p(t) = E_b(t) =\frac{1}{2}\int_0^\ell\{ \vert u_x(x,t)\vert^2+ \vert u_t(x,t) \vert^2 \}dx.
\end{equation}
The main result of this paper is the following.
\begin{theorem} \label{princ}
\begin{enumerate}
\item
We suppose that $\xi = \frac{\ell}{2}$. Then for any $a \in (0, 2)$
there exist  positive constants $C_1, C_2$ such that for all initial data in ${\cal H}$, the  solution of problem
\eqref{a1}-\eqref{a4} satisfies
\begin{equation}\label{expestimate}
E_p(t)\le C_1 e^{- \, C_2t}.
\end{equation}
The constant $C_1$ depends  on the initial data, on $\ell$ and on $a$, while $C_2$
depends only  on $\ell$ and on $a$. 
\item
For any $\mu_1,\mu_2$ satisfying one of the following conditions
\begin{equation}\label{condmu1mu2equiv}
\begin{tabular}{ll}
$1<\mu_2<\mu_1,$\\
$\mu_1<\mu_2<1$,
\end{tabular}
\end{equation}
there exist  positive constants $C_1, C_2$ such that for all initial data in ${\cal H}$, the   solution of problem
\eqref{1.1}--\eqref{1.4} satisfies
\begin{equation}\label{expestimateb}
E_b(t)\le C_1 e^{- \, C_2t}.
\end{equation}
The constant $C_1$ depends  on the initial data, on $\ell$ and on $\mu_1, \mu_2$,
while $C_2$
depends only  on $\ell$ and on $\mu_1, \mu_2$. 
\end{enumerate}
\end{theorem}

The paper is organized as follows. The second section deals with the
well-posedness of the problem while,  in the third section, we prove
the exponential stability of the  systems \rfb{a1}--\rfb{a4} and of \eqref{1.1}--\eqref{1.4}
by using a suitable D'Alembert formula. In  section \ref{section5} we give the same type of results 
for a muldimensional system. 
Some comments and related questions are given in the last section.

\section{Proof of Proposition \ref{propexistunic}}
\setcounter{equation}{0}

Consider the evolution problems
 \be 
\label{OPEN1}
\ddot{y}^j(t) +
 Ay^j(t)  \, 
=  B_1v^j(t), \,\hbox{ in } (2j\ell,2(j+1) \ell), \, j \in \nline^*,
\ee
\be
\label{eq3}
y^j(2j\ell)= \dot{y}^j(2j\ell) = 0, \, j \in \nline^*.
\ee

\be 
\label{eq4}
\ddot{\phi}(t) +
 A\phi(t)  \, 
=  0, \, \hbox{ in }  (0,+\infty),
\ee
\be
\label{OPEN2}
\phi(0)= \phi_0, \, \dot{\phi}(0) = \phi_1.
\ee
A natural question is the regularity of $y^j$ when $v^j \in L^2(2\ell j,2(j+1)\ell j), \, j \in \nline^*$. By 
applying standard energy estimates we can easily check that 
$y^j \in C([2j\ell,2(j+1)\ell];H) \cap C^1([2j\ell,2(j+1)\ell];H_{- \half})$. However if $B_1$ satisfies a 
certain admissibility condition then $y^j$ is more regular. More precisely the 
following result, which is a version of the general transposition method 
(see, for instance, Lions and Magenes \cite{linsmag}) holds true.

It is clear that the system \rfb{eq4}--\rfb{OPEN2} admits a unique solution $\phi$ 
having the regularity 
$$
\phi \in C([0,2 \ell];H_\half) \cap C^{1}([0 j,2  \ell];H),
$$
$$
(\phi,\dot{\phi})(t) = {\mathcal  T}(t)\left(
\begin{array}{ccc}\phi_0\\\phi_1\end{array}\right), \quad
0\leq t \leq 2  \ell.
$$
Moreover,
$B_1^* \phi(\cdot)\in H^1(0,2 \ell), $ and for all $T \in (0, 2\ell)$ there exists a constant
$C>0$
such that 
\be
\Vert (B_1^* \phi)^{\prime}(\cdot)\Vert_{L^2(0,T)}\le C \, \Vert (\phi_0,\phi_1)\Vert_{
H_\half \times H},
\FORALL (\phi_0,\phi_1)\in H_\half \times H.
\label{CACHE1}\ee

\begin{lemma} \label{ex}
Suppose that $v^j \in L^2([2j\ell,2(j+1)\ell]), \, j \in \nline^*$.  
Then the problem \rfb{OPEN1}--\rfb{eq3} admits a unique solution 
having the regularity
\be
y^j \in C([2j\ell,2 (j + 1) \ell];H_\half) \cap C^{1}([2\ell j,2 (j + 1) \ell];H), \  j \in \nline^*,
\label{REG1}\ee
and
$$
(y^{j},\dot{y}^{j})(t) =   \ds \int_{2j\ell}^{t} {\mathcal  T}_{-1}(t-2j\ell-s)\left(
\begin{array}{ccc}0\\B_1 v_j(s)\end{array}\right) \, ds, \quad
2j \ell \leq t \leq 2 (j+1) \ell,\  j \geq 1.
$$
\end{lemma}
\begin{proof}
If we set $Z(t) = \left(
\begin{array}{ll}
y^j(t+2j\ell) \\
\nm
\ds
\dot{y}^j(t+2j\ell)
\end{array}
\right)$ it is clear that \rfb{OPEN1}--\rfb{eq3} can be written as
\[
\dot{Z}^j + {\cal A} Z^j(t) = {\cal B}_1 v^j(t+2j\ell) \hbox{ on } (0, 2\ell), \, Z^j(0) = 0,
\]
where 
\[
{\cal A} = \left(
\begin{array}{cc}
0 & - I \\
\nm
\ds
A &0
\end{array}
\right): H_\half \times H \rightarrow [{\cal D}({\cal A})]^{\prime} ,
\]
\[
{\cal B}_1 = \left(
\begin{array}{ll}
0 \\
\nm
\ds
B_1
\end{array}
\right): \rline \rightarrow [{\cal D}({\cal A})]^{\prime}
\]
It is well known that ${\cal A}$ is a skew adjoint operator so it generates a group of isometries in $[{\cal D}({\cal A})]^{\prime} $, denoted by ${\cal S}(t) (= {\cal T}_{-1}(t))$.

After simple calculations we get that the operator ${\cal B}_1^* : {\cal D}({\cal A}) \rightarrow \rline$ is given by 
\[
{\cal B}_1^*  \left(\begin{array}{c}
u^j \\
v^j
\end{array}
\right) = B^*_1 v^j , \,\forall \, 
(u^j,v^j) \in  {\cal D}({\cal A}).
\]
This implies that 
\[
{\cal B}_1^*{\cal S}^*(t) 
\begin{pmatrix} \phi_0 \cr \phi_1 \end{pmatrix} = B_1^* \dot{\phi}(t), \, \forall \, 
(\phi_0,
\phi_1) \in  {\cal D}({\cal A}),
\]
with $\phi$ satisfying \rfb{eq4}--\rfb{OPEN2}. From the inequality above and 
\rfb{CACHE1} we deduce that there exists a constant $C > 0$ such that 
for all $T\in (0, 2\ell)$
\[
\int_{0}^{T}  \left| {\cal B}_1^*{\cal S}^*(t) \left(\begin{array}{c}
\phi_0 \\
\phi_1
\end{array}
\right) \right|^2 \, dt \le C \, ||(\phi_0,\phi_1)||^2_{H_\half \times H}, \, \forall \, 
(\phi_0,\phi_1) \in  {\cal D}({\cal A}).
\]
According to Theorem 3.1 in \cite[p.187]{ben} (see also \cite{tucsnakweiss}) the inequality above implies 
the interior regularity  \rfb{REG1}.
\end{proof}

The existence result for problem \rfb{a1}--\rfb{a4} is now made by   induction.
First on $[0,2\ell]$ (case   $j=0$), we take
$$
\left(
\begin{array}{ccc}
u^0(t)\\ u^0_t(t)
\end{array}\right) = {\mathcal  T}(t)\left(
\begin{array}{ccc}u_0\\u_1\end{array}\right), \, \forall t \in [0,2\ell].$$
That is clearly a solution of \rfb{a1}--\rfb{a4} on $(0,2\ell)$
  and  that has the regularity $(u^0,u^0_t) \in C([0,2\ell];{\cal H})$.
Now for $j\geq 1$,   we take for all $t\in [2j \ell ,2 (j+1) \ell]$,
\beqs
\left(
\begin{array}{ccc}
u^{j}(t)\\ u_t^{j}(t)
\end{array}\right) &=&\left(
\begin{array}{ccc}
\phi(t+2j\ell)\\ \dot{\phi}(t+2j\ell)
\end{array}\right) + 
\left(
\begin{array}{ccc}
y^{j}(t)\\ \dot{y}^j(t)
\end{array}\right)
\\
&=&{\mathcal  T}(t+2j\ell)\left(
\begin{array}{ccc}u^{j-1}(2j\ell) \\ \dot{u}^{j-1}(2j\ell)\end{array}\right)   +\ds \int_{2j\ell}^{t} {\mathcal  T}_{-1}(t-s)\left(
\begin{array}{ccc}0\\- a \, u_t^{j-1}(s - 2 \ell) \delta_\xi\end{array}\right) \, ds, 
\eeqs
where $y^j$ (resp. $\phi$)  is solution of \rfb{OPEN1}--\rfb{eq3} (resp. \rfb{eq4}--\rfb{OPEN2})
with   $v^j(t) = - a \, u^{j-1}_t(t - 2 \ell)$ (that belongs to $L^2(2j\ell,2(j+1)\ell)$) and   $\phi_0 = u^{j-1}(2j\ell),$  $\phi_1 = \dot{u}^{j-1}(2j\ell)$.
This solution has the announced regularity due to the above arguments.
 
By the same way we prove the second assertion of Proposition \ref{propexistunic}.

\section{Proof of Theorem \ref{princ} }\label{mainproof}
\setcounter{equation}{0}


\noindent
We show that the system \rfb{a1}--\rfb{a4} can be reformulated as follows:
\begin{eqnarray}
& &u^0_{tt}(x,t)- u^0_{xx}(x,t)   = 0 , \ \quad \mbox{\rm in} \ (0,\ell) \times (0,2 \ell),\label{a1bc}\\
& &u^{j-}_{tt}(x,t)- u^{j-}_{xx}(x,t) = 0 ,  \quad \mbox{\rm in }\  
(0,\xi) \times (2 j\ell,2(j+1)\ell),\ j\in \nline^*, \label{a1bb}\\
& &u^{j+}_{tt}(x,t)- u^{j+}_{xx}(x,t) = 0 ,  \quad \mbox{\rm in } \ 
(\xi,\ell) \times (2 j\ell,2(j+1)\ell),\ j\in \nline^*, \label{a1bbb}\\
& &u^{j-}(\xi,t) = u^{j+}(\xi,t),\quad \quad\quad\mbox{\rm on }\  
 (2 j\ell,2(j+1)\ell),\ j \in \nline^*, \label{troppo}\\
& &  - u_x^{j-}(\xi,t) + u_x^{j+}(\xi,t) = - \, a \, u_t^{j+}(\xi,t-2 \ell)\quad \mbox{\rm on }\  
 (2 j\ell,2(j+1)\ell),\ j \in \nline^*,\label{a1bbbc}\\
& &u^{j-}(0,t)=0, \ u^{j+}_x(\ell,t) = 0,\quad \mbox{\rm on}\ \quad(2j\ell,2(j+1) \ell), \, j \in \nline^*, \label{a2bb}\\
& &u(x,0) = u_0(x),  \quad  u_t(x,0)=u_1(x), \quad \mbox{\rm in} \  (0,\ell).\label{a4bb}
\end{eqnarray}

Note that if $u_t^{j+}(\xi,t-2 \ell)$ is replaced by $u_t^{j+}(\xi,t)$, then the energy is decaying if $a> 0$
(and is exponentially decaying if $\xi$ satisfies some conditions: $\xi$ admits a coprime factorization $\frac{p}{q}$ and $p$ is odd).

Hence we look for $u$ solution of (\ref{a1})--(\ref{a4}) in the
form:
\begin{equation}\label{defu-b}
u(x,t)=\alpha_-(x+t)-\alpha_-(t-x),\  \forall x\in (0,\xi), \  t\ge 0,
\end{equation}
and
\begin{equation}\label{defu-}
u(x,t)=\alpha_+(x-\ell+t)+\alpha_+(t-x+\ell),\  \forall x\in (\xi,\ell),\  t\ge 0,
\end{equation}
where $\alpha_-$ and $\alpha_+$ have to be determined.
From this expression we directly see that
$$
u(0,t)=0, \hbox{ and } u_x(\ell,t)=0,\  \forall t\geq 0,
$$
in order words   (\ref{a2}) holds.
Hence it remains to impose the initial conditions at $t=0$ and the transmission conditions at $x=\xi$.

In order to fulfil the
initial conditions (\ref{a4}) for $x\leq \xi$, we take 

\beqs \alpha_-(x)&=&-\frac12
u_0(-x)+\frac12\int_0^{-x} u_1(s)ds \quad \forall
x\in (-\xi,0),\\
\alpha_-(x)&=&\frac12 u_0(x)+\frac12\int_0^x u_1(s)ds \quad \forall
x\in [0,\xi). \eeqs 
In that way $\alpha_-$ is uniquely determined in $(-\xi,\xi)$.

In the same manner   to fulfil the
initial conditions (\ref{a4}) for $x\geq \xi$, we take 

\beqs \alpha_+(y)&=&\frac12
u_0(\ell+y)+\frac12\int_0^{\ell+y} u_1(s)ds \quad \forall
x\in (-(\ell-\xi),0),\\
\alpha_+(y)&=&\frac12 u_0(\ell-y)-\frac12\int_0^{\ell-y} u_1(s)ds \quad \forall
y\in [0,\ell-\xi). \eeqs 
In that way $\alpha_+$ is uniquely determined in $(-(\ell-\xi),\ell-\xi)$.

To check (\ref{a1bc}), we need the continuity of $u$ and $u_x$ at $\xi$, that is equivalent to
\begin{eqnarray*}
\alpha_-(\xi+t)-\alpha_-(t-\xi)=\alpha_+(\xi-\ell+t)+\alpha_+(t-\xi+\ell), \forall t\in (0,2\ell),\\
\alpha'_-(\xi+t)+\alpha'_-(t-\xi)=\alpha'_+(\xi-\ell+t)-\alpha'_+(t-\xi+\ell), \forall t\in (0,2\ell).
\end{eqnarray*}
By setting $y=\xi+t$, this is equivalent to
\begin{eqnarray*}
\alpha_-(y)-\alpha_+(y-2\xi+\ell)=
\alpha_-(y-2\xi)+\alpha_+(y-\ell), \forall y\in (\xi,\xi+2\ell),\\
\alpha'_-(y)+\alpha'_+(y+\ell-2\xi)=-\alpha'_-(y-2\xi)+\alpha'_+(y-\ell), \forall y\in (\xi,\xi+2\ell).
\end{eqnarray*}
Differentiating the first identity in $y$, taking the sum and the difference of the two identities, we get
\begin{eqnarray}\label{serge1}
\alpha'_-(y)=\alpha'_+(y-\ell), \forall y\in (\xi,\xi+2\ell),\\
\alpha'_+(y+\ell-2\xi)=-\alpha'_-(y-2\xi), \forall y\in (\xi,\xi+2\ell).
\label{serge2}
\end{eqnarray}

By iteration this allows to find $\alpha_-$ (resp. $\alpha_+$) up to $2\ell+\xi$ (resp. $3\ell-\xi$).
Indeed fix $\varepsilon\leq 2\min\{\xi,\ell-\xi\}$, then 
in a first step for $y\in (\xi,\xi+\varepsilon)$, we remark that 
$y-\ell$ belongs to $(\xi-\ell,\xi+\varepsilon-\ell)$
which is included in $(\xi-\ell, \ell-\xi)$ the set where $\alpha_+$ is defined up to now. This allows to obtain $\alpha_-^\prime(y)$ for all $y\in (\xi, \xi +\varepsilon ).$
In the same manner $\alpha'_-(y-2\xi)$ is well-defined and this allows then to
obtain $\alpha'_+(y+\ell-2\xi)$ for all $y\in (\xi,\xi+\varepsilon)$.
We now iterate this argument, namely for
$y\in (\xi+\varepsilon,\xi+2\varepsilon)$, the right-hand sides of (\ref{serge1})--(\ref{serge2})
are meaningful, and consequently we obtain 
$\alpha'_-(y)$ (resp. $\alpha'_+(y+\ell-2\xi)$) for such $y$.
We iterate this procedure up to $y\in (\xi+(k-1)\varepsilon,\xi+k\varepsilon)$,
with $k\in \nline$ such that
$$
\xi+k\varepsilon=\xi+2\ell.
$$
This proves the announced statement.

For $y>\xi+2\ell$, we need to take into account (\ref{troppo}) and (\ref{a1bbbc}), that take the form 
$$
\alpha_-(\xi+t)-\alpha_-(t-\xi)=\alpha_+(\xi-\ell+t)+\alpha_+(t-\xi+\ell), \forall t>2\ell,
$$
$$
\alpha'_-(\xi+t)+\alpha'_-(t-\xi)=\alpha'_+(\xi-\ell+t)-\alpha'_+(t-\xi+\ell)
+a(\alpha'_+(\xi+t-3\ell)+\alpha'_+(t-\xi-\ell)), \forall t>2\ell.
$$
By setting $y=\xi+t$, this is equivalent to
$$
\alpha_-(y)-\alpha_+(y-2\xi+\ell)=
\alpha_-(y-2\xi)+\alpha_+(y-\ell), \forall y>\xi+2\ell,
$$
$$
\alpha'_-(y)+\alpha'_+(y+\ell-2\xi)=-\alpha'_-(y-2\xi)+\alpha'_+(y-\ell)
+a(\alpha'_+(y-3\ell)+\alpha'_+(y-2\xi-\ell)), \forall y>\xi+2\ell.
$$
As before differentiating the first equation in $y$ and taking the sum and the difference,
we arrive at (compare with (\ref{serge1})--(\ref{serge2}))
\be
\label{serge3}
\alpha'_-(y)=\alpha'_+(y-\ell)
+\frac{a}{2}(\alpha'_+(y-3\ell)+\alpha'_+(y-2\xi-\ell)), \forall y>\xi+2\ell,
\ee
\be
\label{serge4}
\alpha'_+(y+\ell-2\xi)=-\alpha'_-(y-2\xi)
+\frac{a}{2}(\alpha'_+(y-3\ell)+\alpha'_+(y-2\xi-\ell)), \forall y>\xi+2\ell.
\ee
The same iterative argument allows to show that 
$\alpha_-(y)$ (resp. $\alpha_+(y)$) is uniquely defined for  $y>2\ell+\xi$ (resp. $y>3\ell-\xi$).
Note that this construction based on the D'Alembert formula re-proves the existence result
from Proposition
\ref{propexistunic}. This construction is only valid in one dimension and for a constant coefficients operator, while the semigroup approach of Proposition
\ref{propexistunic} is valid in a more general setting (see below).

The main point is this last iterative relation between $\alpha'_-(y)$, $\alpha'_+(y+\ell-2\xi)$
and previous evaluations.

Let us now take $\xi=\frac{\ell}{2}$, then we can equivalently write 
(\ref{serge3})--(\ref{serge4}) as the following system
\begin{equation}\label{serge5}
\left(
\begin{array}{lll}
\alpha'_-(y)\\
\alpha'_+(y)\\
\alpha'_+(y-\ell)
\\
\alpha'_+(y-2\ell)
\end{array}
\right)=
\left(
\begin{array}{llll}
0&1&\frac{a}{2}&\frac{a}{2}\\
-1&0&\frac{a}{2}&\frac{a}{2}\\
0&1&0&0\\
0&0&1&0
\end{array}
\right)
\left(
\begin{array}{lll}
\alpha'_-(y-\ell)\\
\alpha'_+(y-\ell)\\
\alpha'_+(y-2\ell)
\\
\alpha'_+(y-3\ell)
\end{array}
\right).
\end{equation}
As in \cite{g,gt} we are reduced to calculate the eigenvalues of the matrix
$$
M_a=\left(
\begin{array}{llll}
0&1&\frac{a}{2}&\frac{a}{2}\\
-1&0&\frac{a}{2}&\frac{a}{2}\\
0&1&0&0\\
0&0&1&0
\end{array}
\right)
$$
whose characteristic polynomial is given by
$$
p_a(\lambda)=\lambda^4+(1-\frac{a}{2})\lambda^2+\frac{a}{2}.
$$
The zeroes of $p_a$ are given by
$$
\lambda^2=\frac{a-2\pm\sqrt{a^2-12 a+4}}{4}.
$$
Consequently the eigenvalues of $M_a$ are strictly less than 1 in modulus  
if and only if
\begin{equation}\label{cdspect}
|a-2\pm\sqrt{a^2-12 a+4}|<4.
\end{equation}
In the case $a^2-12 a+4\geq 0$ we see that (\ref{cdspect}) holds
if and only if 
\begin{equation}\label{cdspect1}
0<a\leq 6-4\sqrt{2}.
\end{equation}
On the contrary in the case  $a^2-12 a+4 < 0$ we check that (\ref{cdspect}) holds
if and only if 
\begin{equation}\label{cdspect2}
  6-4\sqrt{2} \leq a<2.
\end{equation}
Hence we conclude that (\ref{cdspect}) holds if and only if $a\in (0,2)$.

Since
$$
p_a'(\lambda)=\lambda (4\lambda^2+2-a),
$$
we can conclude that 
for  $a\in (0, 2)$,  all eigenvalues of $M_a$ are  of modulus $<1$ 
and simple. In that case, there exists a matrix $V_a$ such that
$$
M_a=V_a^{-1} D_a V_a,
$$
where $D_a$ is the diagonal matrix made of the eigenvalues of $M_a$.

Now coming back to (\ref{serge5}) and using an inductive argument, we can deduce that for all $j\in \mathbb N$, 
and for all $y \in (\frac{5\ell}{2}+j\ell, \frac{5\ell}{2}+(j+1)\ell]$, we have
$$
C(y)=M_a^j C(y-j\ell),
$$
where for shortness we have written
$$
C(y):=\left(
\begin{array}{lll}
\alpha'_-(y)\\
\alpha'_+(y)\\
\alpha'_+(y-\ell)
\\
\alpha'_+(y-2\ell)
\end{array}
\right).
$$
Therefore using the previous factorization of $M_a$, we get
$$
C(y)=V_a^{-1} D_a^j V_a C(y-j\ell).
$$
Fnally we find a positive constant $C_a$ (depending only on $a$) such that
for all $j\in \mathbb N$, 
and   all $y \in (\frac{5\ell}{2}+j\ell, \frac{5\ell}{2}+(j+1)\ell]$, we have
\begin{equation}\label{normeit}
 \|C(y)\|_2\leq C_a\rho_a^j \|C(y-j\ell)\|_2,
\end{equation}
where $\rho_a$ is the spectral radius of $D_a$ that is $<1$ (if $a\in (0,2)$).

By simple calculation we see that 
$$
E(t)=\int_{-\frac{\ell}{2}}^{\frac{\ell}{2}} (\alpha'_-(x+t)^2 +\alpha'_+(x+t)^2) \,dx.
$$
Now we closely follow   the arguments of \cite{g,gt}  to conclude the exponential decay of the system.
Namely for all $j\in \mathbb N$, 
and for all $t \in (2\ell+j\ell, 2\ell+(j+2)\ell]$, we can apply (\ref{normeit}) with $y=x+t$
for any $x\in (-\frac{\ell}{2},\frac{\ell}{2})$
and consequently
\begin{eqnarray*}
E(t)&\leq& \int_{-\frac{\ell}{2}}^{\frac{\ell}{2}} \|C(x+t)\|_2^2 \,dx\\
&\leq& C_a^2\rho_a^{2j}  \int_{-\frac{\ell}{2}}^{\frac{\ell}{2}} \|C(x+t-j\ell)\|_2^2 \,dx.
\end{eqnarray*}
Finally as for $t \in (2\ell+j\ell, 2\ell+(j+2)\ell]$ and 
 $x\in (-\frac{\ell}{2},\frac{\ell}{2})$, $x+t-j\ell$ belongs to a compact set,
the quantity
$$
\int_{-\frac{\ell}{2}}^{\frac{\ell}{2}} \|C(x+t-j\ell)\|_2^2 \,dx
$$
is bounded independently of $j$. This means that we have found
a constant $K_a$ such that
for all $j\in \mathbb N$, 
and   all $t \in (2\ell+j\ell, 2\ell+(j+2)\ell]$, one has
$$ 
E(t) \leq  K_a \rho_a^{2j}.
$$ 
This leads to the conclusion because
$\rho_a^{2j}=e^{2j\ln \rho_a}\leq e^{2t\frac{\ln \rho_a}{\ell}}$.

\noindent
Now we study problem (\ref{1.1})--(\ref{1.4})  and  look for  a solution $u$  in the
form:
\begin{equation}\label{defu}
u(x,t)=\alpha(x+t)-\alpha(t-x), \forall x\in (0,\ell),\  t\geq 0.
\end{equation}
Hence we see that (\ref{1.2}) always holds. In order to fulfil the
initial conditions (\ref{1.4}), we take 

\beqs \alpha(x)&=&-\frac12
u_0(-x)+\frac12\int_0^{-x} u_1(s)ds \quad \forall
x\in (-\ell,0),\\
\alpha(x)&=&\frac12 u_0(x)+\frac12\int_0^x u_1(s)ds \quad \forall
x\in [0,\ell). \eeqs To check (\ref{1.3}) we need that
$$
\alpha'(\ell+t)+\alpha'(t-\ell)=0, \hbox{ for } 0<t<2\ell,
$$
or equivalently
$$
\alpha'(y)=-\alpha'(y-2\ell ) \, \forall \, y \in (\ell, 3\ell).
$$
Since the right-hand side is known we get the existence of $\alpha$
on $(\ell, 3\ell)$.

The condition (\ref{dampstand})  is satisfied if
$$
\alpha'(\ell+t)+\alpha'(t-\ell)=\mu_1 (\alpha'(\ell+t)-\alpha'(t-\ell)), \hbox{
for } t\in ((2i+1)2\ell,(2i+2)2\ell),
$$
that is equivalent to
$$
(1-\mu_1)\alpha'(y) =-(1+\mu_1)\alpha'(y-2\ell ), \forall y \in
((2i+1)2\ell+\ell,(2i+2)2\ell+\ell).
$$
Hence for $\mu_1\ne 1$, we find that \be\label{recurrence}
\alpha'(y) =\kappa \alpha'(y-2\ell ), \forall y \in
((2i+1)2\ell+\ell,(2i+2)2\ell+\ell), \ee where $\kappa=\frac{1+\mu_1}{\mu_1-1}$.

In the same manner to check (\ref{dampdelay}) we require
that
$$
\alpha'(\ell+t)+\alpha'(t-\ell)=\mu_2
(\alpha'(t-\ell)-\alpha'(t-3\ell)), \hbox{ for } t\in
((2i+2)2\ell,(2i+3)2\ell),
$$
or equivalently \be\label{DL} \alpha'(y) =(\mu_2-1)\alpha'(y-2\ell)-\mu_2\alpha'(y-4\ell), \forall y \in
((2i+2)2\ell+\ell,(2i+3)2\ell+\ell). \ee

By recurrence we can show that $\alpha$ is well-defined on the whole
$(-\ell, \infty)$.

Now combining  (\ref{recurrence}) and (\ref{DL}) we see that
for $y \in
((2i+1)2\ell+\ell,(2i+2)2\ell+\ell)$ with $i\geq 1$ we obtain
\begin{equation}\label{correct1}
\alpha'(y) =\kappa \alpha'(y-2\ell )=\kappa((\mu_2-1)\alpha'(y-4\ell )-\mu_2\alpha'
(y-6\ell)).
\end{equation}

For $y>7\ell$  we can define the vector
$$U(y):=\left(
\begin{array}{l}
\alpha'(y)\\
\alpha'(y-2l)
\end{array}
\right )
$$
and then,
from  (\ref{DL}) and (\ref{correct1}) we deduce
$$
U(y)=M U(y-4\ell),
$$
where $M$ is the matrix
$$
M=\left(
\begin{array}{ll}
\kappa (\mu_2-1)&-\kappa \mu_2\\
(\mu_2-1)&-  \mu_2
\end{array}
\right).
$$

The eigenvalues of $M$ are $\lambda_1=0$ and $\lambda_2=\kappa(\mu_2-1)-\mu_2.$
Therefore, exponential stability holds if

\begin{equation}\label{condmu1mu2}
\vert \kappa (\mu_2-1)-\mu_2\vert <1.
\end{equation}
Indeed the energy $E_b$  of our system  defined by
$$
E_b(t)=\frac12\int_0^\ell (u_t(x,t)^2+u_x(x,t)^2)\,dx
$$
is here equal to
$$
E_b(t)= \int_{-\ell}^\ell \alpha'(x+t)^2 \,dx.
$$
Hence, the previous arguments show that the energy is exponentially decaying if 
condition (\ref{condmu1mu2}) is satisfied.

Finally by distinguishing the case $\mu_1>1$ to the case $\mu_1<1$, we easily check that
(\ref{condmu1mu2}) is equivalent to (\ref{condmu1mu2equiv}).

\section{The multidimensional case} \label{section5}
\setcounter{equation}{0}
We study the  following internal stabilization problem of a switching delay wave equation in $\Omega\subset \mathbb{R}^d, d\geq 1$. For  given times $T^*>0$ and $\tau\in (0,T^*]$, consider the problem
\be
u_{tt}(x,t) - \Delta u(x,t)+b_1 u_t(x,t)=0\quad \mbox{\rm in}\quad \Omega \times
(i(T^*+\tau), i(T^*+\tau)+T^*),\label{dampstandmultid} \ee 
\be
u_{tt}(x,t) - \Delta u(x,t)+b_2 u_t(x,t-\tau)=0\quad \mbox{\rm in} \quad \Omega \times
( i(T^*+\tau)+T^*, (i+1)(T^*+\tau)),\label{dampdelaymultid} \ee \be
u(x,t)  =0\quad \mbox{\rm on}\quad
\partial \Omega\times (0,\infty)\quad  \label{Dirbc}\ee \be
u(x,0)=u_0(x)\quad \mbox{\rm and}\quad u_t(x,0)=u_1(x)\quad
\hbox{\rm in} \quad \Omega,\label{ic} \ee
where $i\in \nline, b_1>0$ and $b_2$ is a real number.

Note that in the interval $(0,T^*)$ the damping is a standard one, in the sense that it induces an exponential decay of the energy. Hence by standard technique (see e.g. \cite{zuazua,at}), if $T^*$ is fixed such that
the observability estimate in $\Omega$ is valid, there exists $\alpha\in (0,1)$ such that
\begin{equation}\label{decaystandard}
E(T^*)\leq \alpha E(0),
\end{equation}
where $E(t)$ is the standard energy, $E(t):=\frac 1 2\int_{\Omega}(u_t^2+\vert \nabla u\vert^2) dx.$

Now for $t\in (T^*,T^*+\tau)$, by integration by parts, we see that
$$
E'(t)=-\int_\Omega b_2 u_t(x,t)u_t(x,t-\tau)\, dx.
$$
Hence by Cauchy-Schwarz's inequality we find that
$$
E'(t)\leq 2|b_2| E(t)^{1/2}  E(t-\tau)^{1/2}.
$$
Since $t-\tau$ belongs to $(0,T^*)$ and since the energy is decaying on this interval, by (\ref{decaystandard}), we find that
$$
E'(t)\leq 2 \sqrt{\alpha} |b_2| E(0)^{1/2} E(t)^{1/2}.
$$
This can be equivalently written as
$$
\frac{d}{dt } {E(t)}^{1/2}\leq   \sqrt{\alpha} |b_2| E(0)^{1/2},
$$
and integrating this estimate between $T^*$ and $t\in (T^*,T^*+\tau)$, we obtain
$$
{E(t)}^{1/2}-{E(T^*)}^{1/2}\leq \sqrt{\alpha} |b_2| E(0)^{1/2} (t-T^*)\leq \sqrt{\alpha} |b_2| E(0)^{1/2} \tau.
$$
Using again (\ref{decaystandard}), we arrive at
$$
{E(t)}^{1/2} \leq \sqrt{\alpha} (1+|b_2|\tau) E(0)^{1/2}.
$$
As a consequence if the factor  $\tilde \alpha^{1/2}:= \sqrt{\alpha} (1+|b_2|\tau)$ is strictly less than 1,
then we will get a property like (\ref{decaystandard}) but in the interval $(0,T^*+\tau)$, namely
\begin{equation}\label{decayapresdelay}
E(t)\leq \tilde \alpha E(0), \forall t\in (T^*,T^*+\tau).
\end{equation}
Note that the condition $\tilde \alpha^{1/2}<1$ is equivalent to
\begin{equation}\label{condb2}
|b_2|<\frac{1-\sqrt{\alpha}}{\sqrt{\alpha} \tau},
\end{equation}
that means that
$b_2$ has to be small enough.
Since our system is invariant by a translation of $T^*+\tau$, this argument may be repeated between
$i(T^*+\tau)$ and  $(i+1)(T^*+\tau)$, and therefore we find that
$$
E(t)\leq \tilde \alpha^{i+1} E(0), \forall t\in (i(T^*+\tau),(i+1)(T^*+\tau)).
$$
Writing $\tilde \alpha^{i+1}=e^{(i+1)(T^*+\tau)\frac{\log \tilde \alpha}{(T^*+\tau)}}$
and using the fact that $\frac{\log \tilde \alpha}{(T^*+\tau)}<0$, we arrive at
\begin{equation}\label{decayapresdelayfinal}
E(t)\leq e^{t\frac{\log \tilde \alpha}{(T^*+\tau)}} E(0), \forall t\in (i(T^*+\tau),(i+1)(T^*+\tau)),
\end{equation}
which proves the exponential decay of the energy.
In conclusion we have proved the next result.

\begin{theorem} \label{princ2}
Assume that $T^*$ is the minimal time of observability for the wave equation with internal damping,
that $\tau\in (0,T^*]$
and that $(\ref{condb2})$ holds.
Then the energy of the system $(\ref{dampstandmultid})-(\ref{ic})$ decays exponentially to zero.
\end{theorem}

\br{\rm
1. Our arguments also hold if we replace the internal damping in $\Omega \times
(2iT^*, (2i+1)T^*)$ by a boundary damping. 
Similarly the global internal damping can be replaced by a local one, as far as the exponential decay is guaranteed. Obviously in both cases, the time $T^*$ of observability has to be changed.
The converse situation, namely keep internal damping in $\Omega \times
(2iT^*, (2i+1)T^*)$ and take a boundary damping with delay in $ \Omega \times
((2i+1)T^*, 2(i+1)T^*)$ is more delicate because we are not able to prove that
$$
\int_{\partial \Omega} b_2 u_t(x,t)u_t(x,t-\tau)\, dx\leq 2 |b_2| E(t)^{1/2}  E(t-\tau)^{1/2}.
$$
Hence another argument should be found.

2. Instead of taking a constant coefficient $b_2$ we can also take   $b_2\in L^\infty(\Omega).$ In this case
the condition (\ref{condb2}) has to be replaced by
$$
\sup_\Omega|b_2|<\frac{1-\sqrt{\alpha}}{\sqrt{\alpha} \tau}.
$$
}
\er


\section{Comments and related questions} 
\setcounter{equation}{0}
\begin{enumerate}
\item
The statement of Theorem \ref{princ} concerning problem \rfb{a1}--\rfb{a4}
remains valid in the case $\frac{\xi}{\ell}=\frac{p}{q}$ with 
$p,q\in \mathbb N$ with $p$ odd and $q$ even.
We did not give its proof since it is too technical
and do not bring any new ideas.
We have chosen $\xi = \frac{\ell}{2}$ because this is the best 
location for the decay rate in the absence of delay.
\item
In the same manner we can obtain the same result as Theorem \ref{princ} for the following problem:
\be
u_{tt}(x,t)- u_{xx}(x,t)   = \mu_1 u_t(\xi,t) \delta_\xi ,\  \mbox{\rm in} \quad(0,\ell) \times 
(2i \ell,2(i+1)\ell),\forall i\in \nline,\label{a1switch} \ee $$
u_{tt}(x,t)- u_{xx}(x,t) + a \, u_t(\xi,t-2 \ell) \, \delta_\xi = 0,\hspace{4.7cm} $$
\be
\hspace{5 cm}  \mbox{\rm in } (0,\ell) \times (2(i+1)\ell, 2(i+2)\ell),\forall i\in \nline,\label{a1bswitch} \ee \be
u(0,t)=0, \  u_x(\ell,t) = 0,\  \mbox{\rm on}\quad(0,+\infty), \label{a2switch} \ee \be
u(x,0) = u_0(x),  \quad  u_t(x,0)=u_1(x), \  \mbox{\rm in} \quad (0,\ell).\label{a4switch}
\ee

\item
Let $H$ be a Hilbert space equipped with the norm $||.||_{H}$, and
let $A :\cD(A)\rightarrow H$ be a self-adjoint, positive and
invertible operator. We introduce the scale of Hilbert spaces
$H_{\alpha}$, $\alpha\in\rline$, as follows\m: for every $\alpha\geq
0$, $H_{\alpha}=\Dscr(A^{\alpha})$ with the norm $\|z
\|_{\alpha}=\|A^\alpha z\|_{H}$. The space $H_{-\alpha}$ is defined
by duality with respect to the pivot space $H$ as  $H_{-\alpha}
=H_{\alpha}^*$ for $\alpha>0$. The operator $A$ can be extended (or
restricted) to each $H_{\alpha}$ such that it becomes a bounded
operator \BEQ{A0ext} A : H_{\alpha} \rarrow H_{\alpha - 1} \,
\,\mbox{for}\,\, \alpha\in\rline \m. \EEQ

The second ingredient needed for our construction is a bounded
linear operator $B~ :~ U \longrightarrow ~H_{-\half},$ where $U$ is
another Hilbert space identified with its dual. The operator $B^*$
is bounded from $H_{\frac12}$ to $U$.

The systems that  we considered in this paper enter in one of the following abstract problems:
\BEQ{damped1} \ddot w(t)
+ Aw(t)  
= 0, \quad 0 \leq t\leq  T_0, \EEQ 
\BEQ{int}
\ddot w(t)
+ Aw(t)  + \mu \, BB^*\dot w(t- T_0) = 0, \quad t\geq  T_0,
\EEQ
\BEQ{output} w(0) \m=\m w_0,  \m \dot
w(0)=w_1, \EEQ 
or
\BEQ{damped1b} \ddot w(t)
+ Aw(t)  
= 0, \quad 0 \leq t\leq  T_0, \EEQ 
\BEQ{damped1bb} \ddot w(t)
+ Aw(t) + \mu_1 BB^* \dot{u}(t)
= 0, \quad  (2i + 1)T_0 \leq t\leq (2i + 2)T_0, \, 
\forall \, i \in \nline, \EEQ 
\BEQ{intb}
\ddot w(t)
+ Aw(t)  + \mu_2 \, BB^*\dot w(t-T_0) = 0, \quad (2i + 2)T_0 \leq t\leq (2i + 3)T_0, \, \forall \, i \in \nline,
\EEQ
\BEQ{outputb} w(0) \m=\m w_0,  \m \dot
w(0)=w_1, \EEQ
where
$T_0
> 0$ is the time delay, $\mu, \mu_1, \mu_2$ are real numbers and the initial datum $(w_0,w_1)$ belongs to a suitable
space.

Assume that there exist $T \geq T_0,C > 0$ such that 
\be 
\label{inegobs}
\int_{0}^{T}  \left\|B^{*} \phi^\prime(s)\right\|_{U}^2 d\,s \leq C\,
||(w_0,w_1)||^2_{H_\half \times H}
\ee
for $(w_0,w_1) \in H_1 \times H_\half$ and $\phi$ is the solution of
the undamped evolution equation
\be \ddot \phi(t) + A \phi(t) =
0, \quad t\geq0, \label{eq3b} \ee \be \phi(0) = w_0,
 \dot \phi(0) = w_1.\nonumber
  \label{eq4b}\ee
 
To study the well--posedness of the system
\eqref{damped1}--\eqref{output}, we write it as an abstract
Cauchy problem in a product Banach space, and use the semigroup
approach. For this take the Hilbert space ${\mathcal H}~:=~
H_{\half} \times H$ and the unbounded
linear operators
 \be
  \label{opbb} {\mathcal  A}: {\mathcal
D}({\mathcal  A}) = H_1 \times H_\half \subset {\mathcal  H} \longrightarrow {\mathcal  H}, \, {\mathcal
A}\left(
\begin{array}{ccc}u_1\\u_2\end{array}\right) = \left(
\begin{array}{l}
u_2 \\
- Au_1  
\end{array}
\right) \ee

and

 $$ {\mathcal  A}_d: {\mathcal
D}({\mathcal  A}_d) = \left\{(u,v) \in {\cal H}; \, v \in H_\half, \, Au + \mu_1 BB^*v \in H \right\} \subset {\mathcal  H} \longrightarrow {\mathcal  H}, \, 
$$
\be
  \label{opbbb}
{\mathcal
A}_d \left(
\begin{array}{ccc}u_1\\u_2\end{array}\right) = \left(
\begin{array}{l}
u_2 \\
- Au_1  - \mu_1 \, BB^*u_2
\end{array}
\right). \ee

The operators $({\mathcal  A}, {\mathcal  D}({\mathcal  A}) )$ and $({\mathcal  A}_d, {\mathcal  D}({\mathcal  A}_d) )$  defined by \eqref{opbbb} 
generate a strongly continuous semigroup of contractions on ${\mathcal H}$ denoted respectively by  $({\mathcal
T}(t))_{t\geq0}$ and $({\mathcal
T}_d(t))_{t\geq0}$ (as before let $({\mathcal
T}_{-1}(t))_{t\geq0}$ be the extension of $({\mathcal
T}(t))_{t\geq0}$ to $H_{-1}$).
 
\begin{proposition}\label{propexistunicb}
\begin{enumerate}
\item
Assume that the inequality \rfb{inegobs}   holds. Then the system
\rfb{damped1}--\rfb{output} is well--posed. More precisely, for every
$(u_0,u_1)\in {\mathcal  H}$, the solution of \rfb{damped1}--\rfb{output} is given by 
$$
\left(
\begin{array}{ccc}
u(t)\\ \dot{u}(t) 
\end{array}\right)= \left\{
\begin{array}{ll}
\left(
\begin{array}{ccc}
u^0(t)\\ \dot{u}^0(t)
\end{array}\right) = {\mathcal  T}(t)\left(
\begin{array}{ccc}u_0\\u_1\end{array}\right), \, 0 \leq t \leq  T_0, \\
\left(
\begin{array}{ccc}
u^{j}(t)\\ \dot{u}^{j}(t)
\end{array}\right) = {\mathcal  T}(t-j T_0)\left(
\begin{array}{ccc}u^{j-1}(j T_0)\\ \dot{u}^{j-1}(j T_0)\end{array}\right) + \\ \ds \int_{jT_0}^{t} {\mathcal  T}_{-1}(t-s)\left(
\begin{array}{ccc}0\\ - \mu \, BB^* \dot{u}^{j-1}(s)\end{array}\right) \, ds, \, \\
\hspace{5cm} j T_0 \leq t \leq (j+1) T_0, j \geq 1.

\end{array}
\right.
$$ 
and satisfies $(u^j,\dot{u}^j) \in C([jT_0, (j+1)T_0], {\mathcal  H}), \, j \in \nline.$

\item
Assume that the inequality \rfb{inegobs} holds. Then, the system
\rfb{1.1b}--\rfb{1.4b} is well--posed. More precisely, for every
$(u_0,u_1)\in {\mathcal  H}$, the solution of \rfb{1.1b}--\rfb{1.4b} is given by 
$$
\left(
\begin{array}{ccc}
u(t)\\ u_t(t)
\end{array}\right) = \left\{
\begin{array}{ll}
\left(
\begin{array}{ccc}
u^0(t)\\ u_t^0(t)
\end{array}\right) = {\mathcal  T}(t)\left(
\begin{array}{ccc}u_0\\u_1\end{array}\right), \, 0 \leq t \leq T_0, \\
\left(
\begin{array}{ccc}
u^{2j+1}(t)\\ u_t^{2j+1})(t)
\end{array}\right) = {\mathcal  T}_d(t-(2j+1)T_0)\left(
\begin{array}{ccc}u^{2j}((2j+1)T_0)\\ u_t^{2j}((2j+1)T_0)\end{array}\right), \, 
\\
\hspace{5cm} (2j+1)T_0 \leq t \leq  (2j +2)T_0, j \in \nline,\\
\left(
\begin{array}{ccc}
u^{2j+2}(t)\\ u_t^{2j+2}(t)
\end{array}\right) = {\mathcal  T}(t-(2j+2)T_0)\left(
\begin{array}{ccc}
u^{2j+1}((2j+2)T_0)\\ u_t^{2j+1}((2j+2)T_0)
\end{array}\right) + \\
\ds \int_{2(2j+2) \ell}^{t} {\mathcal  T}_{-1}(t-s)\left(
\begin{array}{ccc}0\\- \mu_2 BB^*u_t^{2j+1}(s - 2 \ell) \delta_\ell\end{array}\right) \, ds, \, \\
\hspace{5cm} (2j+2)T_0 \leq t \leq   (2j+3)T_0 , j \in \nline
\end{array}
\right.
$$ 
and satisfies $$(u^0,u^0_t) \in C([0, 2 \ell], {\mathcal  H}), \, (u^{2j+1},u^{2j+1}_t) \in C([2(2j+1)\ell, 2(2j+2)\ell], {\mathcal  H}), \,j \in \nline, \, $$ 
$$ (u^{2j+2},u^{2j+2}_t) \in C([2(2j+2)\ell, 2(2j+3)\ell], {\mathcal  H}), \, j \in \nline.$$
\end{enumerate}
\end{proposition}

We now give two multi--dimensional illustrations of this setting.
Let $\Omega \subset\RR^n$ be an open bounded set   with a smooth
boundary $\Gamma$.  We assume that $\Gamma$ is divided into two
parts $\Gamma_0$ and $\Gamma_1$, i.e. $\Gamma
=\Gamma_0\cup\Gamma_1,$ with $\overline \Gamma_0 \cap\overline
\Gamma_1 =\emptyset$ and $meas \Gamma_1\neq 0$ (and satisfied some Lions geometric condition or some geometric control condition, see \cite{ammari} and \cite{ak} for more details). Note that the
condition $\overline \Gamma_0 \cap\overline \Gamma_1 =\emptyset$ is
only made in order to simplify the presentation, hence our  analysis
can be performed without this assumption in a similar manner.

We further fix a time interval $T_0>0$ and a delay $\tau>0$. In this
domain $\Omega$  we consider the initial boundary value problems with
switching boundary conditions:

 \begin{eqnarray}
& &u_{tt}  -\Delta u=0\quad \mbox{\rm in}\quad\Omega\times
(0,+\infty),\label{1.1b}\\
& & u   =0\quad \mbox{\rm on}\quad \Gamma_0\times
(0,+\infty),\quad  \label{1.2b}\\
& & u   =0\quad \mbox{\rm on}\quad \Gamma_1\times
(0,T_0),\quad  \label{1.2bb}\\
 & & u(x,t) =\mu_1\frac{\partial G(u_t)}{\partial n}(x,t)  \quad \mbox{\rm
on}\quad\Gamma_1\times
((2i+1)T_0,(2i+2)T_0),\label{dampstandb}\\
& & u(x,t)= \mu_2\frac{\partial G(u_t)}{\partial n}(x, t-\tau)  \quad
\mbox{\rm on}\quad\Gamma_1\times
((2i+2)T_0,(2i+3)T_0),\label{dampdelayb}\\
& &u(x,0)=u_0(x)\quad \mbox{\rm and}\quad u_t(x,0)=u_1(x)\quad
\hbox{\rm in}\quad\Omega,\label{1.4b}
\end{eqnarray}
and
\begin{eqnarray}
& &u_{tt}  + \Delta^2 u=0\quad \mbox{\rm in}\quad\Omega\times
(0,+\infty),\label{1.1bv}\\
& & u   =0\quad \mbox{\rm on}\quad \partial \Omega \times
(0,+\infty),\quad  \label{1.2bv}\\
& & \Delta u = 0 \, \quad \mbox{\rm on}\quad \Gamma_0\times
(0,+\infty),\quad  \label{1.2bvv}\\
& & \Delta u   =0\quad \mbox{\rm on}\quad \Gamma_1\times
(0,T_0),\quad  \label{1.2bbv}\\
 & & \Delta u(x,t) = - \mu_1\frac{\partial G(u_t)}{\partial n}(x,t)  \quad \mbox{\rm
on}\quad\Gamma_1\times
((2i+1)T_0,(2i+2)T_0),\label{dampstandbv}\\
& & \Delta u(x,t)= - \mu_2\frac{\partial G(u_t)}{\partial n}(x, t-\tau)  \,
\mbox{\rm on}\quad\Gamma_1\times
((2i+2)T_0,(2i+3)T_0),\label{dampdelaybv}\\
& &u(x,0)=u_0(x)\quad \mbox{\rm and}\quad u_t(x,0)=u_1(x)\quad
\hbox{\rm in}\quad\Omega,\label{1.4bv}
\end{eqnarray}
where $G = (-\Delta)^{-1} : H^{-1}(\Omega) \rightarrow H^{1}_0(\Omega), i \in \nline, \mu_1$ and $\mu_2$ are real parameters.

Note that the above systems are exponentially stable in absence of
time delay,  that is if $\tau=0$  and if $\mu_1=\mu_2>0$.

According to \cite{ammari} and \cite{ak} the inequality \rfb{inegobs} is satisfied for some $T_0 >0$ and Proposition \ref{propexistunicb} implies that the wave system (\ref{1.1b})--(\ref{1.4b}) and the plate system (\ref{1.1bv})--(\ref{1.4bv}) admit a finite energy solution.

For $\tau = T_0$, where $T_0$ is fixed such that the observability estimate in $\Gamma_1$ is valid, by the same method used in Theorem \ref{princ2}, we can prove an exponential stability result   for both systems.

\end{enumerate}

\end{document}